\documentclass[11pt]{amsart}
\usepackage{graphicx}

\usepackage{amsmath,amsfonts,amsthm,amssymb,graphics,color, hyperref}

\newtheorem{theorem}{Theorem}

\newtheorem{ex}{Exercise}
\newtheorem{question}{Question}

\theoremstyle{remark}

\theoremstyle{definition}

\numberwithin{equation}{section}

\newcommand{\fp}{\mathfrak{P}}

\newcommand{\fS}{\mathfrak{S}} 
\newcommand{\R}{\mathbb{R}}

\renewcommand{\epsilon}{\varepsilon}

        \definecolor{pink}{rgb}{1,0,1}

\title[Love Games]{Love Games \\ A Game Theory Approach to Compatibility}  
\author{Kerstin Bever} \address{Georg-August-Universit\"at G\"ottingen \\ Mathematisches Institut \\ Bunsenstr. 3--5 \\ 37073 G\"ottingen} \author{Julie Rowlett} \address{Leibniz Universit\"at Hannover \\ Institut f\"ur Analysis \\ am Welfengarten 1 \\ 30167 Hannover \\  rowlett@math.uni-hannover.de} 

\begin{document}
\maketitle
\begin{abstract}
In this note, we present a compatibility test based on John Nash's game-theoretic notion of equilibrium strategy.  The test must be taken separately by both partners, making it difficult for either partner alone to control the outcome.  The mathematics behind the test including Nash's celebrated theorem and an example from the film, ``A Beautiful Mind,'' are discussed as well as how to customize the test for more accurate results and how to modify the test to evaluate interpersonal relationships in other settings, not only romantic.  To investigate the long-term dynamics of give and take in a relationship we introduce the ``iterated dating dilemma'' and apply the notion of ``zero-determinant payoff strategy" introduced by Dyson and Press in 2012 for the iterated prisoner's dilemma.  
\end{abstract} 

\section{Introduction}

``Are you ready to settle down?''  ``Are you in love or forcing it?''  These are the titles of compatibility quizzes found in Cosmopolitan magazine \cite{cosmo}.   The ``Love Calculator'' featured in Glamour magazine \cite{glamour} claims to determine whether ``you and he add up" or are ``destined for a long division?''  These and most other compatibility quizzes we have seen in popular culture do not appear to have a rigorous scientific basis and are often easy to cheat to control the test results.  

We have created a compatibility test based on sound mathematical principles and John Nash's notion of equilibrium strategy.  The test must be taken separately by both partners, making it difficult for either partner alone to control the outcome.  

This note is organized as follows.  In \S 2, we introduce Nash's notion of equilibrium strategy and investigate an example from popular culture.  Our test comprises \S 3.  In \S 4, we show how to evaluate the test, and we also show how the test can be modified to provide more accurate results or to evaluate non-romantic relationships, for example in professional or political settings.   In the last section of this note, we explore the long-term dynamics of give and take in a relationship using a variation of the iterated prisoner's dilemma (IPD) and Dyson's and Press's 2012 discovery of ``zero determinant payoff strategies for the IPD \cite{dp}. 

\section{The mathematics behind the test}  
Incorporating rigorous pure mathematics into real-world situations can be a subtle problem.  To illustrate this let's consider an example from popular culture.  

\subsection{A Beautiful Mind} 
Our quiz is based on the notion of equilibrium strategy and the Nobel-prize winning theorem of John Nash \cite{nash}.  To familiarize the uninitiated with this notion, we recall a scene from the film, ``A Beautiful Mind,'' based on Nash's life and work, which attempts to depict Nash's thereom.  The scene is set in a bar and begins with Nash's character and his buddies watching a group of attractive women enter the bar.  One woman is exceptionally beautiful, and the men begin discussing their best strategies to court her.  After a moment of reflection Nash's character realizes that the ``best'' tactic is for the men to each approach one of the averagely attractive women, rather than trying to court the exceptionally beautiful woman.  One sees this depicted in the film as each man leaves with one of the averagely beautiful women, and the exceptionally beautiful woman is alone.  Is this really an equilibrium strategy according to Nash's Theorem?

Let's recall the definitions and statement of Nash's celebrated theorem.  An $n$-person non-cooperative game is based on the absence of coalitions or communication between players.  The $i^{th}$ player has a set of $m_i$ pure strategies, each of which can be identified with one of the standard unit vectors in $\R^{m_i}$.  The set of mixed strategies corresponds to a probability distribution over pure strategies and can be identified with the convex subset 
$$\fS_i \cong \left\{ x = (x_1, \ldots, x_{m_i}) \in \R^{m_i} : x_j \geq 0 \forall j, \quad \sum_{j=1} ^{m_i} x_j =1 \right\}.$$
The total strategy space over all players can then be identified with 
$$\fS \cong \prod_{i=1} ^n \fS_i \subset \R^{N}, \quad N = \sum_{i=1} ^n m_i.$$
Each player has an associated payoff function $\wp_i : \fS \to \R$.  The payoff functions are linear in the strategy of the respective player.  That is, if all other players' strategies are fixed, the payoff function is a linear function from $\fS_i \to \R$.  Typically, the payoff function is assumed to be continuous on $\fS$.  In \cite{nash}, Nash proved that for every such game, there exists at least one ``equilibrium strategy,'' which in the following sense is the ``best'' strategy for all players.  For $s \in \fS$, let $\wp_i (s; \sigma; i)$ denote the payoff for the strategy in which the $i^{th}$ player's strategy according to $s$ is replaced by $\sigma \in \fS_i$.  An equilibrium strategy $s$ satisfies 
$$\wp_i (s) \geq \wp_i (s; \sigma_i; i) \quad \forall \sigma_i \in \fS_i, \quad \forall i=1, \ldots, n.$$
In other words, if a single player changes his strategy, he cannot increase his payoff.  

\begin{theorem}[Nash \cite{nash}] For any $n$-person non-cooperative game such that the payoff functions are linear in the strategy of each player and are continuous functions on the total strategy space there exists at least one equilibrium strategy. 
\end{theorem} 

To analyze the scene from the film, for simplicity, let's assume there are 2 men, denoted by man 1 and man 2, and 3 women, one denoted by ``G'' (for gorgeous) and two denoted by ``P'' (for pretty).  Each man has two pure strategies:  ``G'' which corresponds to courting the gorgeous woman, and ``P'' which corresponds to courting one of the pretty women.  First, consider the following payoff matrix:

\bigskip 
\begin{center} 
\begin{tabular}{|c|c|c|} 
 \hline & G & P \\ \hline G & (-1, -1) & (1, 0) \\ \hline P & (0, 1) & (0, 0) \\ \hline 
\end{tabular}
\end{center} 
 \bigskip 

The payoff matrix means:  
$$\wp_1 (G, G) = \wp_2 (G, G) = -1, \wp_1 (G, P) = 1, \wp_2 (G, P) = 0,$$
$$ \wp_1 (P, G) = 0, \wp_2 (P, G) = 1, \wp_1 (P, P) = \wp_2 (P, P) =0.$$
The interpretation is that if both men attempt to court the gorgeous woman, they are both unsuccessful, indicated by a negative payoff $-1$.   If man 1 chooses ``G,'' and man 2 chooses ``P,'' then both men's courtships are successful, and we consider man 1 to be positively rewarded indicated by a score of 1, and man 2 is merely content, indicated by a neutral score of 0.  If both men choose ``P'' we presume that they each court a different pretty woman and are both successful, so they are both content, indicated by a score of 0.  The game is symmetric.  It's straightforward to generalize to $n+1$ women and $n$ men, but this is not necessary to illustrate our point.  

By Nash's Theorem, there exists at least one equilibrium strategy, and in this simple example it's possible to compute the equilibrium strategy or strategies.  If the probabilities that man 1 and 2 choose ``G'' are denoted by $x$, and $y$, respectively, then the total strategy space is 
$$\fS \cong \{ (x,y) \in \R^2 : 0 \leq x,y \leq 1 \}.$$
The payoff functions are 
$$\wp_1 (x,y) = x - 2xy, \quad \wp_2 (x,y) = y - 2xy.$$

\begin{ex} Determine the equilibrium strategies.   
\end{ex}

Hopefully we've all arrived at the same answer to the exercise:  the equilibrium strategies are $(x=0, y=1)$ and $(x=1, y=0)$.  Does this fit with the film?  The first equilibrium strategy means that man 1 courts a pretty woman, and man 2 courts the gorgeous woman.  The second equilibrium strategy is the other way around.  This doesn't fit with the film at all!  According to the film, the equilibrium strategy should be that both men court pretty women, which corresponds to $x=y=0$. 

This is rather perplexing.  One explanation is that Hollywood just doesn't understand math.  That may very well be the case.  However, as mathematicians, we ought to know that things aren't always what they seem.  Is it possible to resolve the mathematical definition of equilibrium strategy with the scene depicted in the film?    

\begin{ex} What is a payoff matrix for a two-player symmetric game with two pure strategies G and P with the interpretations described above, such that the equilibrium strategy is $x=y=0$?  
\end{ex}

How do we begin?  Let's denote the pure strategies as above by G and P.  The payoffs are symmetric, so  $\wp_1 (G, P) = \wp_2 (P, G)$, $\wp_1 (P, P) = \wp_2 (P, P)$, and $\wp_1 (G, G) = \wp_2 (G, G)$.  The unknown payoff matrix for the film is: 

\bigskip 
\begin{center} 
\begin{tabular}{|c|c|c|} 
 \hline & G & P \\  \hline G & (a, a) & (c, d) \\ \hline P & (d, c) & (b, b) \\ \hline 
\end{tabular}
\end{center} 
 \bigskip 

The payoff functions are 
$$\wp_1 (x,y) = xya + x(1-y)c + (1-x)yd + (1-x)(1-y)b,$$
$$\wp_2 (x,y) = xya + y(1-x)c + (1-y)xd + (1-y)(1-x)b.$$
The payoffs for the strategy $x=y=0$ are $\wp_1 (0,0) = b = \wp_2 (0,0)$.  In order for this to be an equilibrium strategy, the following must hold: 
$$\wp_1 (x,0) \leq \wp_1 (0,0), \quad \forall x \in [0,1]; \quad \wp_2 (0,y) \leq \wp_2 (0,0), \quad \forall y \in [0,1].$$
This is equivalent to: 
$$xc + (1-x) b \leq b \quad \forall x \in [0,1], \quad yc + (1-y)b \leq b \quad , \forall y \in [0,1],$$
which holds iff 
$$c \leq b.$$
This means that the payoff associated to successfully courting the gorgeous woman is \em less than or equal to \em the payoff associated to successfully courting a pretty woman!  This would seem inconsistent with the film which indicates that each man would be most content if he successfully courted the gorgeous woman.

This leads us to the same initial conclusion:  Hollywood simply doesn't understand the math.  Well, perhaps there is a logical explanation.  What if $b=c$? The payoff matrix would be 

\bigskip 
\begin{center}
\begin{tabular}{|c|c|c|} 
 \hline & G & P \\  \hline G & (a, a) & (c, d) \\ \hline P & (d, c) & (c, c) \\ \hline 
\end{tabular}
\end{center} 
 \bigskip 

Now, clearly $a<c$, because successfully courting the gorgeous woman is a more positive outcome than unsuccessfully courting her.  It is also natural to presume that $d<c$, as the film indicates that successfully courting the gorgeous woman while the other man courts a pretty woman is more desirable than the other way around.  How can we make sense of the fact that $\wp_1(P, P) = \wp_2(P, P) = \wp_1 (G, P) = \wp_2 (P, G) = c$?  Instead of thinking about the relationship between the men and the women, let's think about the relationship \em between the men.  \em  They are shown as buddies in the film.  If man 1 courts the gorgeous woman while his buddy, man 2, courts a pretty woman, then although man 1 might be very happy about the outcome, man 2 might be jealous.  Jealously can have unpleasant consequences.  Now, on the other hand, if both men court a pretty woman, then nobody will be jealous!  So, perhaps jealously explains why the above payoff matrix with $d<a<c$ is what the filmmakers had in mind.  There's still a glitch with this explanation.   You'll see what we mean if you do the following exercise.

\begin{ex} Are there any \em other \em equilibrium strategies?  
\end{ex}

Well, yes.  The strategy $(G, G)$ is also an equilibrium strategy, and this is clearly not the best outcome since both men end up alone.  So, it's still not entirely clear whether or not the filmmakers understand the mathematics.  Now let's move on to our test which \em is \em consistent with rigorous mathematics!   
\section{The Compatibility Test}
\subsection{Instructions}
Instead of simply deciding on one of three options as is typical for love quizzes, both partners have to  
fill in a matrix consisting of three slots:   

\begin{minipage}[t]{0.2\textwidth}
\begin{tabular}[c]{c |p{0.7 cm}| } 
\multicolumn{2}{r}{Partner \# 1} \\ \cline{2-2}
\(a\) & \hspace{2cm}\\ \cline{2-2}
\(b\) & \hspace{2cm}\\  \cline{2-2}
\(c\) & \hspace{2cm}\\  \cline{2-2}
\end{tabular}
\end{minipage}
\begin{minipage}[t]{0.2\textwidth}
\begin{tabular}[h]{c |p{0.7 cm}| } 
\multicolumn{2}{r}{Partner \# 2} \\ \cline{2-2}
\(\alpha\) & \hspace{2cm}\\ \cline{2-2}
\(\beta\) & \hspace{2cm}\\  \cline{2-2}
\(\gamma\) & \hspace{2cm}\\  \cline{2-2}
\end{tabular}
\end{minipage}
\begin{minipage}[t]{0.1\textwidth}
\begin{tabular}[h]{l} \vspace{0.1cm}\\
Outcome A\\ 
Outcome B \\
Outcome C \\ 
\end{tabular}
\end{minipage}\\

Partners \# 1 and \#2 each \em separately \em choose a number between 0 and 10 to fill in each of the slots \((a/\alpha,b/ \beta, c / \gamma )\), such that the numbers add up to 10.  For example, you could choose 0, 1, 9; or 2, 3, 5; or 3, 3, 4.  How do you choose the numbers?  

\begin{enumerate}
\item For the $a/\alpha$ slot, choose a number from 0 (never/completely disagree) to 10 (always/completely agree) to describe how likely Outcome A is.   
\item For the $b/\beta$ slot, choose a number between 0  (never/completely disagree) to 10 (always/completely agree) to describe how likely Outcome B is.   
\item For the $c/\gamma$ slot, choose a number between 0  (never/completely disagree) to 10 (always/completely agree) to describe how likely Outcome C is.    
\end{enumerate}

It is important when taking the test that the partners do not see each other's answers!  For this reason, we recommend printing out two copies of the following section and separately filling in your numbers, leaving the numbers of your partner's responses blank.  

\subsection{The test questions}

\begin{question} 
You have decided to spend the weekend together. How does it go?  \\

\begin{minipage}[t]{0.2\textwidth}
\begin{tabular}[c]{c |p{0.7 cm}| } 
\multicolumn{2}{r}{Partner \# 1} \\ \cline{2-2}
\(a\) & \hspace{2cm}\\ \cline{2-2}
\(b\) & \hspace{2cm}\\  \cline{2-2}
\(c\) & \hspace{2cm}\\  \cline{2-2}
\end{tabular}
\end{minipage}
\begin{minipage}[t]{0.2\textwidth}
\begin{tabular}[h]{c |p{0.7 cm}| } 
\multicolumn{2}{r}{Partner \# 2} \\ \cline{2-2}
\(\alpha\) & \hspace{2cm}\\ \cline{2-2}
\(\beta\) & \hspace{2cm}\\  \cline{2-2}
\(\gamma\) & \hspace{2cm}\\  \cline{2-2}
\end{tabular}
\end{minipage}
\begin{minipage}[t]{0.1\textwidth}
\begin{tabular}[h]{l} \vspace{0.1cm}\\
It's a blast!\\ 
It's not my first choice, but it's nice.\\
Well, at least my partner is happy.\\
\end{tabular}
\end{minipage}\\
\end{question} 

\begin{question} 
How do you feel about your sex life?  \\

\begin{minipage}[t]{0.2\textwidth}
\begin{tabular}[c]{c |p{0.7 cm}| } 
\multicolumn{2}{r}{Partner \# 1} \\ \cline{2-2}
\(a\) & \hspace{2cm}\\ \cline{2-2}
\(b\) & \hspace{2cm}\\  \cline{2-2}
\(c\) & \hspace{2cm}\\  \cline{2-2}
\end{tabular}
\end{minipage}
\begin{minipage}[t]{0.2\textwidth}
\begin{tabular}[h]{c |p{0.7 cm}| } 
\multicolumn{2}{r}{Partner \# 2} \\ \cline{2-2}
\(\alpha\) & \hspace{2cm}\\ \cline{2-2}
\(\beta\) & \hspace{2cm}\\  \cline{2-2}
\(\gamma\) & \hspace{2cm}\\  \cline{2-2}
\end{tabular}
\end{minipage}
\begin{minipage}[t]{0.1\textwidth}
\begin{tabular}[h]{l} \vspace{0.1cm}\\
Sex is just the way I like it! \\ 
I'm satisfied.\\
I am sexually frustrated.\\
\end{tabular}
\end{minipage}\\
\end{question} 

\begin{question} 
How do you and your partner manage your careers and chores?   \\

\begin{minipage}[t]{0.2\textwidth}
\begin{tabular}[c]{c |p{0.7 cm}| } 
\multicolumn{2}{r}{Partner \# 1} \\ \cline{2-2}
\(a\) & \hspace{2cm}\\ \cline{2-2}
\(b\) & \hspace{2cm}\\  \cline{2-2}
\(c\) & \hspace{2cm}\\  \cline{2-2}
\end{tabular}
\end{minipage}
\begin{minipage}[t]{0.2\textwidth}
\begin{tabular}[h]{c |p{0.7 cm}| } 
\multicolumn{2}{r}{Partner \# 2} \\ \cline{2-2}
\(\alpha\) & \hspace{2cm}\\ \cline{2-2}
\(\beta\) & \hspace{2cm}\\  \cline{2-2}
\(\gamma\) & \hspace{2cm}\\  \cline{2-2}
\end{tabular}
\end{minipage}
\begin{minipage}[t]{0.1\textwidth}
\begin{tabular}[h]{l} \vspace{0.1cm}\\
Career and chores are just the way I like them.\\  
I've had to make some compromises.  \\
I've made significant sacrifices for my partner.\\
\end{tabular}
\end{minipage}\\
\end{question} 

\begin{question} 
Your partner has fallen ill.  What happens?  \\

\begin{minipage}[t]{0.2\textwidth}
\begin{tabular}[c]{c |p{0.7 cm}| } 
\multicolumn{2}{r}{Partner \# 1} \\ \cline{2-2}
\(a\) & \hspace{2cm}\\ \cline{2-2}
\(b\) & \hspace{2cm}\\  \cline{2-2}
\(c\) & \hspace{2cm}\\  \cline{2-2}
\end{tabular}
\end{minipage}
\begin{minipage}[t]{0.2\textwidth}
\begin{tabular}[h]{c |p{0.7 cm}| } 
\multicolumn{2}{r}{Partner \# 2} \\ \cline{2-2}
\(\alpha\) & \hspace{2cm}\\ \cline{2-2}
\(\beta\) & \hspace{2cm}\\  \cline{2-2}
\(\gamma\) & \hspace{2cm}\\  \cline{2-2}
\end{tabular}
\end{minipage}
\begin{minipage}[t]{0.1\textwidth}
\begin{tabular}[h]{l} \vspace{0.1cm}\\

I rarely catch whatever s/he has. \\ 
I take care of him/her and might get sick.\\
I stay by his/her side and get sick. \\
\end{tabular}
\end{minipage} \\
\end{question} 

\begin{question} 
How often do you see your family and your in-laws? \\

\begin{minipage}[t]{0.2\textwidth}
\begin{tabular}[c]{c |p{0.7 cm}| } 
\multicolumn{2}{r}{Partner \# 1} \\ \cline{2-2}
\(a\) & \hspace{2cm}\\ \cline{2-2}
\(b\) & \hspace{2cm}\\  \cline{2-2}
\(c\) & \hspace{2cm}\\  \cline{2-2}
\end{tabular}
\end{minipage}
\begin{minipage}[t]{0.2\textwidth}
\begin{tabular}[h]{c |p{0.7 cm}| } 
\multicolumn{2}{r}{Partner \# 2} \\ \cline{2-2}
\(\alpha\) & \hspace{2cm}\\ \cline{2-2}
\(\beta\) & \hspace{2cm}\\  \cline{2-2}
\(\gamma\) & \hspace{2cm}\\  \cline{2-2}
\end{tabular}
\end{minipage}
\begin{minipage}[t]{0.1\textwidth}
\begin{tabular}[h]{l} \vspace{0.1cm}\\
Exactly as much as I would like.  \\
It's a compromise but we manage to get along.\\
I don't see my own family enough. \\ 
\end{tabular}
\end{minipage} \\
\end{question} 
 
\newpage
\begin{question} 
If you are in a long relationship, your circle of friends often change.  Whose friends do you tend to see?  \\
 
\begin{minipage}[t]{0.2\textwidth}
\begin{tabular}[c]{c |p{0.7 cm}| } 
\multicolumn{2}{r}{Partner \# 1} \\ \cline{2-2}
\(a\) & \hspace{2cm}\\ \cline{2-2}
\(b\) & \hspace{2cm}\\  \cline{2-2}
\(c\) & \hspace{2cm}\\  \cline{2-2}
\end{tabular}
\end{minipage}
\begin{minipage}[t]{0.2\textwidth}
\begin{tabular}[h]{c |p{0.7 cm}| } 
\multicolumn{2}{r}{Partner \# 2} \\ \cline{2-2}
\(\alpha\) & \hspace{2cm}\\ \cline{2-2}
\(\beta\) & \hspace{2cm}\\  \cline{2-2}
\(\gamma\) & \hspace{2cm}\\  \cline{2-2}
\end{tabular}
\end{minipage}
\begin{minipage}[t]{0.1\textwidth}
\begin{tabular}[h]{l} \vspace{0.1cm}\\
I spend as much time with my friends as I want. \\ 
I spend less time with my friends but still keep up.\\
I rarely see my old friends.\\ 
\end{tabular}
\end{minipage} \\
\end{question}

\begin{question} 
How is your financial situation?   \\

\begin{minipage}[t]{0.2\textwidth}
\begin{tabular}[c]{c |p{0.7 cm}| } 
\multicolumn{2}{r}{Partner \# 1} \\ \cline{2-2}
\(a\) & \hspace{2cm}\\ \cline{2-2}
\(b\) & \hspace{2cm}\\  \cline{2-2}
\(c\) & \hspace{2cm}\\  \cline{2-2}
\end{tabular}
\end{minipage}
\begin{minipage}[t]{0.2\textwidth}
\begin{tabular}[h]{c |p{0.7 cm}| } 
\multicolumn{2}{r}{Partner \# 2} \\ \cline{2-2}
\(\alpha\) & \hspace{2cm}\\ \cline{2-2}
\(\beta\) & \hspace{2cm}\\  \cline{2-2}
\(\gamma\) & \hspace{2cm}\\  \cline{2-2}
\end{tabular}
\end{minipage}
\begin{minipage}[t]{0.1\textwidth}
\begin{tabular}[h]{l} \vspace{0.1cm}\\
Great.\\ 
I can't always spend the way I'd like, but it's fine. \\
We have financial problems. \\ 
\end{tabular}
\end{minipage} \\
\end{question} 

\begin{question} 
Have you and your partner talked about having children?  \\

\begin{minipage}[t]{0.2\textwidth}
\begin{tabular}[c]{c |p{0.7 cm}| } 
\multicolumn{2}{r}{Partner \# 1} \\ \cline{2-2}
\(a\) & \hspace{2cm}\\ \cline{2-2}
\(b\) & \hspace{2cm}\\  \cline{2-2}
\(c\) & \hspace{2cm}\\  \cline{2-2}
\end{tabular}
\end{minipage}
\begin{minipage}[t]{0.2\textwidth}
\begin{tabular}[h]{c |p{0.7 cm}| } 
\multicolumn{2}{r}{Partner \# 2} \\ \cline{2-2}
\(\alpha\) & \hspace{2cm}\\ \cline{2-2}
\(\beta\) & \hspace{2cm}\\  \cline{2-2}
\(\gamma\) & \hspace{2cm}\\  \cline{2-2}
\end{tabular}
\end{minipage}
\begin{minipage}[t]{0.1\textwidth}
\begin{tabular}[h]{l} \vspace{0.1cm}\\
Yes and I'm happy with our decisions.\\ 
We will discuss it eventually. \\
That is up to my partner to decide.\\ 
\end{tabular}
\end{minipage} \\
\end{question} 

\begin{question} 
How do you feel about your lifestyle in terms of health, fitness, and physical appearance?   \\

\begin{minipage}[t]{0.2\textwidth}
\begin{tabular}[c]{c |p{0.7 cm}| } 
\multicolumn{2}{r}{Partner \# 1} \\ \cline{2-2}
\(a\) & \hspace{2cm}\\ \cline{2-2}
\(b\) & \hspace{2cm}\\  \cline{2-2}
\(c\) & \hspace{2cm}\\  \cline{2-2}
\end{tabular}
\end{minipage}
\begin{minipage}[t]{0.2\textwidth}
\begin{tabular}[h]{c |p{0.7 cm}| } 
\multicolumn{2}{r}{Partner \# 2} \\ \cline{2-2}
\(\alpha\) & \hspace{2cm}\\ \cline{2-2}
\(\beta\) & \hspace{2cm}\\  \cline{2-2}
\(\gamma\) & \hspace{2cm}\\  \cline{2-2}
\end{tabular}
\end{minipage}
\begin{minipage}[t]{0.1\textwidth}
\begin{tabular}[h]{l} \vspace{0.1cm}\\
I am totally happy. \\ 
Fine. \\
It is not what it used to be.  \\ 
\end{tabular}
\end{minipage} \\
\end{question} 

\begin{question} 
You have agreed to spend a cozy night at home watching TV.  How does it go? \\

\begin{minipage}[t]{0.2\textwidth}
\begin{tabular}[c]{c |p{0.7 cm}| } 
\multicolumn{2}{r}{Partner \# 1} \\ \cline{2-2}
\(a\) & \hspace{2cm}\\ \cline{2-2}
\(b\) & \hspace{2cm}\\  \cline{2-2}
\(c\) & \hspace{2cm}\\  \cline{2-2}
\end{tabular}
\end{minipage}
\begin{minipage}[t]{0.2\textwidth}
\begin{tabular}[h]{c |p{0.7 cm}| } 
\multicolumn{2}{r}{Partner \# 2} \\ \cline{2-2}
\(\alpha\) & \hspace{2cm}\\ \cline{2-2}
\(\beta\) & \hspace{2cm}\\  \cline{2-2}
\(\gamma\) & \hspace{2cm}\\  \cline{2-2}
\end{tabular}
\end{minipage}
\begin{minipage}[t]{0.1\textwidth}
\begin{tabular}[h]{l} \vspace{0.1cm}\\
I love what we watch.  \\  
We compromise on something.  \\
My partner chooses what we watch.\\ 
\end{tabular}
\end{minipage}
\end{question}

\section{How to evaluate the test} 

The notion of equilibrium strategy was introduced for \em non-cooperative games.  \em You may be thinking
\begin{quote} \em Why use \textbf{non}-cooperative game theory for a \textbf{compatibility} test?  \em  
\end{quote} 
Although it may seem counter-intuitive, there is a quite natural reason.  Most everyday decisions are made spontaneously with little or no discussion.  People naturally tend to act in their best self-interest.  This is similar to non-cooperative games in which each player acts independently to maximize his own payoff without communication with the other players.   

The quiz is based on a simple, two-player, two-move, symmetric, zero-sum game with one dominant strategy.  We use this game to quantify a snap decision made by a couple without discussion.  Strategy ``M'' corresponds to ``my way,'' which means that player does exactly what s/he would like to do.  Strategy ``H'' corresponds to ``the highway,'' which means that player does the opposite of what s/he would like to do.  The payoff $\wp_1(M, M) = \wp_2 (M, M)$ means that if both partners make a snap decision which fits with each of their personal interests, then they are both equally content, and neither one nor the other is dominant, expressed by a neutral payoff (0) to both partners.  The payoffs $\wp_1 (M, H) = 1, \wp_2 (M, H) = -1$ means that if the decision favors Partner \# 1's interest and is not consistent with Partner \# 2's interest, Partner \# 1 has a positive payoff +1, whereas Partner \# 2 has a negative payoff -1.  If both partners' interests are not met, then neither is dominant, expressed by a neutral payoff of $0$ to each.  

\begin{figure} 
\begin{center} 
\begin{tabular}{|c|c|c|} 
 \hline & M & H \\ \hline M & (0, 0) & (1, -1) \\ \hline H & (-1, 1) & (0, 0) \\ \hline 
\end{tabular}
 \caption{My way or the highway?}  \label{quizfig}
\end{center} 
\end{figure} 

How can such a simple game quantify happiness and balance in a relationship?  
\begin{ex} What is the equilibrium strategy for the above game?
\end{ex} 

The unique equilibrium strategy is $(M, M)$.  This corresponds to both partners having it ``my way.''  Well, that's clearly the best possible scenario!  This would mean that neither partner dominates, and yet both are totally happy with their life choices and activities together.  We designed the test to quantify the following two characteristics of a relationship.   
\begin{enumerate}
\item Is the relationship balanced, or is one of the partners dominant?  If so, who?  
\item How happy are both partners with their lives in the relationship? 
\end{enumerate}

How does our test relate to the model game?  Let $x$ denote the probability that Partner \# 1 has it ``my way" and $y$ denote the probability that Partner \# 2 has it ``my way."  Then corresponding payoffs to each partner, 
$$\wp_1 (x,y) = x-y, \quad \wp_2 (x,y) = y-x.$$
The test questions were designed to measure the cumulative probability that each partner has it ``my way'' or ``the highway.''  For a general test with $10$ questions on a 10-point scale as in \S 3, we compute 
$$X := \sum_{t=1}^{10} a(t), \quad  \textrm{and} \quad Y:= \sum_{t=1}^{10} \alpha(t).$$ 

What do the above sums mean?  If we divide each by $100$, then $x = X/100$ and $y = Y/100$ are the cumulative probabilities that each partner has it his or her way.  The payoffs for the game ``my way or the highway" are then 
$$\wp_1 = x-y, \quad \wp_2 = y-x.$$
Since it is easier to work with whole numbers, we use $X$ and $Y$ and compute 
$$\fp_1 := X - Y, \quad \fp_2 := Y - X.$$
\begin{enumerate}
\item If $\fp_1 > \fp_1$, then the first partner is dominant.  The larger the difference, the more partner \# 1 dominates.  
\item If $\fp_2 > \fp_1$, then the second partner is dominant.  The larger the difference, the more partner \# 2 dominates
\end{enumerate} 

Your relationship could be in balance, but that doesn't have to mean that you are happy.  For example, if both partners answer $0$ to each $a/\alpha$ slot and answers $10$ to each $c/\gamma$ slot, that means that although neither dominates, neither partner is pleased with the outcome of the situations posed in the test questions.  
This is where the second parameter comes into play.  The more both of you feel that you're doing the things you like to do, the more satisfied you are. This is why it makes sense to look at the following numbers:  
$$ \mathfrak{K}_1:= \sum_{t=1} ^{10} a (t) - c (t), \quad \textrm{and} \quad \mathfrak{K}_2 := \sum_{t=1}^{10} \alpha (t) - \gamma (t).$$

\begin{enumerate}
\item The larger $\mathfrak{K}_1$ is, the more the first partner feels that s/he does what s/he likes. 
\item The larger $\mathfrak{K}_2$ is, the more the second partner feels that s/he does what s/he likes. 
\item The overall satisfaction in the a relationship is given by 
$$\mathfrak{K} := \mathfrak{K}_1 + \mathfrak{K}_2.$$
The larger this number is, the more both partners are doing what they'd like to do. 
\end{enumerate}

\begin{figure} 
\includegraphics[width=160pt]{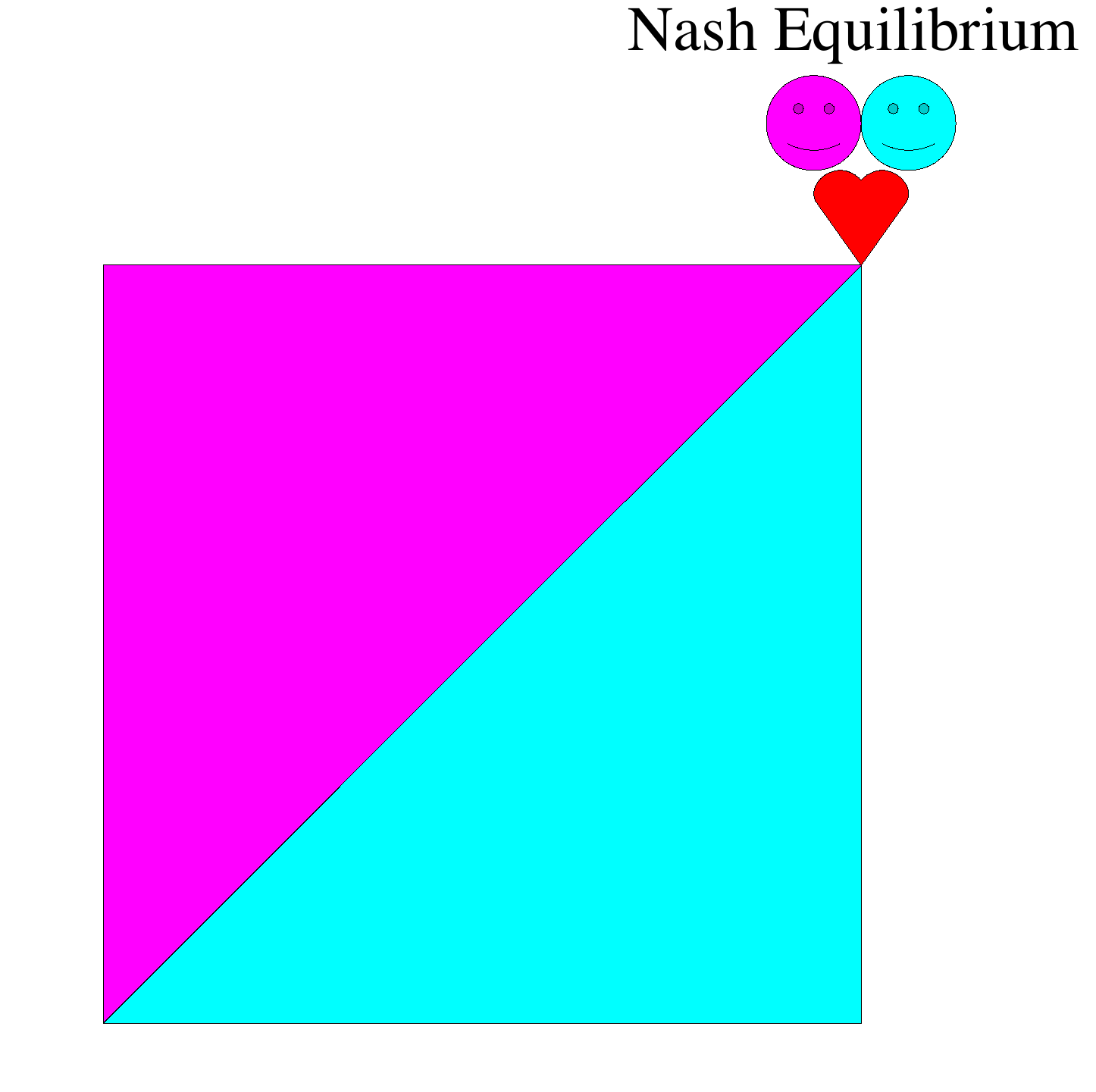} \caption{For the two-player, symmetric, zero-sum game ``my way or the highway,'' the unique equilibrium strategy is when both players execute ``my way'' with probability 1.  Letting the respective probabilities that partners 1 and 2 execute ``my way'' be $(x, y)$, then in the magenta region $y>x$, which corresponds to $\wp_2 > \wp_1$ in the test, and the second partner is dominant.  In the cyan region, $x>y$, which corresponds to $\wp_1 > \wp_2$ in the test, and the first partner is dominant.  The relationship is balanced along the diagonal line $x=y$.  Both partners are maximally happy at the point $(1,1)$, which corresponds to the maximum possible value of $\mathfrak{K}$.} \label{neq-fig} 
\end{figure} 

The unique Nash equilibrium strategy for the model game corresponds to $\fp_1 = \fp_2 = 0$, and $\mathfrak{K}$ having the maximum possible value which for the test presented here is 200.  Our test measures how close a couple is to this perfectly balanced happiness (see Figure \ref{neq-fig}).  

\subsection{Modifying the test}  
The two main parameters quantified by our test are:
\begin{enumerate}
 \item Is the relationship balanced, or is one partner more dominant?  If so, who?  
\item How satisfied are both partners in the relationship?  
\end{enumerate}

We chose questions to assess sex-life, family and friends, health and lifestyle, career, and finances.  For more accurate results, the test can be altered to consist of as many questions on as many topics as desired.  A further variation of the test is to allow each partner to choose the total of the answer column for each question.  The default total is 10, but if a certain question is more important to the partner, s/he can make that question's column have a higher sum like 50.  On the other hand, if a certain question is really not so important to him or her, then s/he could assign the column's total to be a smaller number like 5.  

To evaluate this type of test one first computes 
$$X := \sum_{t=1} ^T a(t), \quad  \textrm{and} \quad Y:= \sum_{t=1} ^T \alpha(t),$$
where $T$ is the total number of test questions.  Next, these are normalized by the weights of the questions.  If Partner \# 1 assigns $p(t)$ points to question $t$, and Partner \# 2 assigns $q(t)$ points to question $t$, we compute
$$W_1 := \sum_{t=1} ^T p(t), \quad W_2 := \sum_{t=1} ^T q(t).$$
Then we define 
$$\wp_1 := \frac{X}{W_1} - \frac{Y}{W_2}, \quad \wp_2 := \frac{Y}{W_2} - \frac{X}{W_1}.$$
The numbers $\frac{X}{W_1}$ and $\frac{Y}{W_2}$ are, respectively, the cumulative probabilities that Partners \# 1 and \# 2 execute strategy ``my way'' in the ``my way or the highway" game, and $\wp_1$, $\wp_2$ are their respective payoffs.  

\begin{enumerate}
\item If $\wp_1 > \wp_2$, then the first partner dominates, and the larger $\wp_1 - \wp_2$, the more Partner \# 1 dominates. 
\item If $\wp_2 > \wp_1$, then the second partner dominates, and the larger $\wp_2 - \wp_1$, the more Partner \# 2 dominates. 
\end{enumerate}

To compute the overall satisfaction of the pair, we compute 
$$ \mathfrak{K}_1 := \sum_{t=1} ^{T} \frac{a (t) - c (t)}{p(t)}, \quad \textrm{and} \quad \mathfrak{K}_2 := \sum_{t=1}^{T} \frac{\alpha (t) - \gamma (t)}{q(t)}.$$
\begin{enumerate}
\item The larger $\mathfrak{K}_1$ is, the more the first partner feels that s/he does what s/he likes. 
\item The larger $\mathfrak{K}_2$ is, the more the second partner feels that s/he does what s/he likes. 
\item The overall satisfaction in the a relationship is given by 
$$\mathfrak{K} := \mathfrak{K}_1 + \mathfrak{K}_2.$$
The larger this number is, the more both partners are doing what they'd like to do based on what they feel is most important.  In this case, the largest possible values of $\mathfrak{K}_1$ and $\mathfrak{K}_2$ are both $T$, and the largest possible value of $\mathfrak{K}$ is therefore $2T$.  
\end{enumerate}

With a suitable choice of test questions, our test could also be used to evaluate relationships between business partners, employees, political partners, or any situation involving two (or more) people who behave in their best self interest.  It is also possible to create tests based on iterated games as in \cite{bev}.

\section{The Iterated Dating Dilemma}  
Have you been in a relationship for a length of time?  In any relationship there is give and take.  We can describe this mathematically using the prisoner's dilemma.  Two crooks have committed a crime and have gotten caught!  They are locked in a prison cell together before they are taken separately for questioning.  They agree in the prison cell that they won't rat each other out.  This is known as \em cooperation (C).  \em  However, when each prisoner goes off for questioning, he has the chance to \em defect (D) \em by claiming he is innocent, and his partner is solely responsible for the crime.  

The payoff matrix (Figure \ref{odd}) means that if both prisoners cooperate, they both receive an equal payoff of $W$.  In terms of give and take within a couple, we identify the strategy C with \em giving, \em and the strategy D with \em taking.  \em  So if both partners give, then they both receive an equal payoff $W$.  In the prisoner dilemma, if one prisoner cooperates but his partner in crime defects, then the prisoner who cooperated gets slapped with a longer prison sentence than the defector.  So the payoff for cooperating $Z$ is smaller than $Y$, the payoff for defecting.  Similarly, if one person in a relationship gives, corresponding to strategy C, while the other person in the relationship takes, corresponding to strategy D, then the payoff to the giver is lower than the payoff to the taker.  Finally, if both prisoners defect, then they both get a longer sentence than if they cooperate, so the payoff for mutual defection $X$ is smaller than the payoff for mutual cooperation $W$.  Similarly, if both people in the relationship take, then their payoff is lower than if they both give.  Consequently, the payoffs satisfy 
$$Y>W>X>Z.$$
In the prisoner's dilemma it is also customary to assume that the expected value of cooperating is higher than the expected value of defecting.  If a prisoner defects, a neutral expectation of his partner in crime is 1/2 probability defecting and 1/2 probability cooperating.  Hence, the expected payoff of defecting is 
$$\frac{1}{2} \left( Y+Z \right),$$
and we assume the payoff for cooperating is better than the expected payoff of defecting so 
$$W > \frac{1}{2} \left( Y + Z \right).$$

\begin{figure} 
\begin{center} 
\begin{tabular}{|c|c|c|} 
 \hline & C & D \\ \hline C & (W,W) & (Z,Y) \\ \hline D & (Y,Z) & (X,X) \\ \hline 
\end{tabular}
\caption{Payoff matrix for the prisoner's dilemma and the IDD}   \label{odd} 
\end{center} 
 \end{figure} 

In the iterated prisoner's dilemma the ``game'' is repeated.  We can use the dynamics of the iterated prisoner's dilemma to study the dynamics of give and take in a relationship, which we call the ``iterated dating dilemma (IDD).''  Of course, we are simplifying things a bit here, for example by assuming the payoffs are symmetric and remain constant over time.  Nonetheless we have found that analyzing the long-term dynamics of give and take in a relationship using the iterated prisoner's dilemma leads to interesting results.  

At each round of the game we assume that each person remembers what happened in the previous round.  Press and Dyson proved that it doesn't matter whether one person can remember more than the other, which game-theoretically means that a long-memory player has no advantage over a memory-one player (someone who can only remember the previous date) \cite{dp}.  This in itself is rather interesting and a bit of a relief!  It means that to ensure our best possible payoff over time, we just need to remember what happened last time.  

Let's call the two people in the couple Pat and Gene.  At each iteration of the IDD there are four possibilities for the previous outcome:  $(C, C)$, $(C, D)$, $(D, C)$, $(D, D)$.  Pat's strategy is $p = (p_1, p_2, p_3, p_4)$ which corresponds to her probability of cooperating (giving) under each of the previous outcomes.  Gene's strategy is $q = (q_1, q_2, q_3, q_4)$ (c.f. Fig. 1 \cite{dp}) which corresponds to his probability of cooperating under each of the previous outcomes.   Hauert and Schuster \cite{hs} showed that these probabilities can be used to determine the expected outcome for each player.  The unexpected fact proven in \cite{dp} is that there exist strategies for Pat and Gene such that Pat can control Gene's expected payoff, and/or Gene can control Pat's expected payoff.  These are known as \em zero-determinant \em (ZD) strategies, because these strategies are precisely those such that the determinant of a certain matrix vanishes; we refer interested readers to the proofs and calculations in \cite{dp} and focus here on the implications of those results.  Press and Dyson showed that Pat (or Gene) can force any particular score for Gene (or Pat) regardless of what the other person does.  They also showed that it is \em not \em possible for Pat or Gene to unilaterally set his or her own score.  

Although it is impossible for Pat to control her own expected payoff, she can control Gene's, and vice-versa.  This would seem to indicate that whereas one can influence his or her partner's happiness by his or her choice of behaviors, one cannot unilaterally control one's own happiness in a relationship.  Moreover, since a longer memory player has no advantage over a memory-one player, it suffices to remember what happened ``last time'' 
to make one's best decision for ``next time.''  

\section*{Acknowledgments} We are grateful to Henry Segerman for critical reading and constructive criticism of this manuscript.  The support of the Georg-August-Universit\"at G\"ottingen is gratefully acknowledged as well as the Max Planck Institut f\"ur Mathematik and the Leibniz Universit\"at Hannover.


\end{document}